\documentclass[12pt]{article}

\usepackage{indentfirst}
\usepackage{amsfonts}
\usepackage{amsmath}
\usepackage{amssymb}
\usepackage{amsbsy, amsthm}

\newtheorem{dfn}{Definition}
\newtheorem{theorem}[dfn]{Theorem}
\newtheorem{lemma}[dfn]{Lemma}

\newenvironment{pf}{\noindent{\bf Proof.}}
{\enspace\vrule height5pt depth0pt width5pt}
\def\deg {{\rm deg}}
\def\lc {{\rm lc}}

\begin{document}
\title{Linear colorings of subcubic graphs}
\author{Chun-Hung Liu$^1$\thanks{E-mail:cliu87@math.gatech.edu.} and  
        Gexin Yu$^2$\thanks{E-mail:gyu@wm.edu. 
		Research supported in part by NSA grant H98230-12-1-0226} \\
 {\small $^1$School of Mathematics, Georgia Institute of Technology, Atlanta, Georgia 30332, USA} \\
 {\small $^2$Department of Mathematics, The College of William and Mary, Williamsburg, VA, 23185, USA}}

\maketitle

\abstract{
A linear coloring of a graph is a proper coloring of the vertices of the graph so that each pair of color classes induce a union of disjoint paths.  
In this paper, we prove that for every connected graph with maximum degree at most three and every assignment of lists of size four to the vertices of the graph, there exists a linear coloring such that the color of each vertex belongs to the list assigned to that vertex and the neighbors of every degree-two vertex receive different colors, unless the graph is $C_5$ or $K_{3,3}$.   
This confirms a conjecture raised by Esperet, Montassier, and Raspaud (L. Esperet, M. Montassier, and A. Raspaud, {\it Linear choosability of graphs},  Discrete Math. 308 (2008) 3938--3950).
Our proof is constructive and yields a linear-time algorithm to find such a coloring.

\bigskip

}

\section{Introduction}

A {\em proper coloring} of a graph is an assignment of colors to the vertices of the graph so that adjacent vertices receive different colors.  Graph coloring is an important topic in graph theory and has wide applications in scheduling and partitioning problems.  

Yuster (1998,~\cite{Y98}) introduced the notion of {\it linear coloring}, which is a proper coloring such that each pair of color classes induce a linear forest, where a linear forest is a union of disjoint paths.  This notion combines the well-studied {\em acyclic colorings} (which are proper colorings so that each pair of color classes induce a forest) introduced by Gr\"{u}nbaum (1973, \cite{G73}) and the frugal colorings (a proper coloring is {\em $k$-frugal} if the subgraph induced by each pair of color classes has maximum degree less than $k$) introduced by Hind, Molloy, and Reed (1997, \cite{HMR97}). 

We write $\lc(G)$ to denote the {\it linear chromatic number} of $G$, which is the smallest integer $k$ such that $G$ has a linear coloring with $k$ colors.  Yuster~\cite{Y98} constructed an infinite family of graphs such that $\lc(G)\ge C_1\Delta(G)^{3/2}$, for some constant $C_1$.  He also proved an upper bound of $\lc(G)\le C_2\Delta(G)^{3/2}$, for some constant $C_2$ and for sufficiently large $\Delta(G)$.

As most coloring problems, it is hard in general to determine the linear chromatic number.  
For example, Esperet, Montassier, and Raspaud~\cite{EMR08} proved that deciding whether a bipartite subcubic graph is linearly $3$-colorable is an NP-complete problem.    
On the other hand, there are some easy upper and lower bounds on $\lc(G)$ for every graph $G$.  
Let $G$ be a graph with maximum degree $\Delta(G)$.  
Then $\lc(G)\ge\lceil{\Delta(G)/2}\rceil+1$, since each color can appear on at most two neighbors of a vertex of maximum degree,  and $\lc(G)\le \chi(G^2)\le\Delta(G^2)+1\le \Delta(G)^2+1$, where $\chi(G)$ denotes the chromatic number of $G$, and $G^2$ is the graph obtained from $G$ by adding edges $xy$ for each pair of vertices $x,y$ with distance two.   
Li, Wang, and Raspaud~\cite{LWR11} improved the easy upper bound to $\lc(G)\le (\Delta(G)^2+\Delta(G))/2$.  

For every family of lists $L = \{L(v): v \in V(G)\}$ of size $k$, we say that a proper coloring $f$ is a {\it proper $L$-coloring} of $G$ if $f(v) \in L(v)$ for every vertex $v$ of $G$.
General list-coloring was first introduced by Erd\"{o}s, Rubin, and Taylor~\cite{ERT79} and independently by Vizing~\cite{V76} in the 1970s, and it has been well-explored since then~\cite{JT95}.  
In particular, the following analog of Brooks' Theorem for list-coloring was proven \cite{ERT79,V76}, and Skulrattanakulchai \cite{S06} gave a linear-time algorithm to find a proper $L$-coloring when the family of lists $\{L(v): v \in V(G)\}$ is given.

\begin{theorem}  \label{list Brooks}
Every connected graph with maximum degree $\Delta(G)$ is $\Delta(G)$-choosable unless $G$ is an odd cycle or a complete graph.
\end{theorem}

The list-version of linear coloring was first studied by Esperet, Montassier, and Raspaud~\cite{EMR08}.  
We say that a proper $L$-coloring of $G$ is {\it linear} if the subgraph of $G$ induced by each two color classes is a linear forest.
A graph $G$ is {\em linearly $k$-choosable} if for every family of lists $(L(v): v \in V(G))$ of size $k$, the graph $G$ has a linear $L$-coloring.
When all the lists are the same, it is the same as linear $k$-coloring.   
We denote by $\lc_l(G)$ the smallest $k$ so that $G$ is linearly $k$-choosable.  

Clearly, $\lc_l(G) \geq \lc(G)$.
Substantial work has been done on the study of graphs whose linear (list) chromatic number is close to the easy lower bound $\lceil \Delta/2 \rceil+1$, see ~\cite{CY11, CH11, DXZ10, EMR08, LWR11, RW09, WW12}.    
On the other hand, a little more is known when a graph has small maximum degree. Li, Wang, and Raspaud~\cite{LWR11} showed that $\lc(G)\le 8$ if $\Delta(G)\le 4$ and $\lc(G)\le 14$ if $\Delta(G)\le 5$.   
Esperet, Montassier, and Raspaud~\cite{EMR08} proved that $\lc_l(G)\le 9$ if $\Delta(G)\le 4$ and $\lc_l(G)\le 5$ if $G$ is subcubic (i.e., $\Delta(G)\le 3$).  
Note that $K_{3,3}$ is not linearly $4$-colorable, so their upper bound on subcubic graph is tight, but they conjectured $K_{3,3}$ is the only subcubic graph which is not linearly $4$-choosable.  
In this paper, we confirm this conjecture. 
As a matter of fact, we prove the following slightly stronger result.  

We say that a linear $L$-coloring is {\it superlinear} if the neighbors of every vertex of degree two receive different colors.
We say that a graph $G$ is {\it superlinearly $k$-choosable} if it is linearly $k$-choosable in such a way that the corresponding linear coloring can be chosen to be superlinear.

\begin{theorem} \label{main}
Let $G$ be a subcubic graph which has no component isomorphic to $K_{3,3}$ or $C_5$. 
Then $G$ is superlinearly $4$-choosable.
\end{theorem}

Note that $C_5$ is linearly $4$-choosable but not superlinearly $4$-choosable.
However, there is a superlinear $L$-coloring of $C_5$ when some vertices $u$ and $v$ have the different lists of colors $L(u)$ and $L(v)$.
In addition, our proof of Theorem \ref{main} is constructive and yields a linear-time algorithm to find a superlinear $L$-coloring when the family of lists $\{L(v): v \in V(G)\}$ is given.
This generalizes an algorithm of Skulrattanakulchai \cite{S04} to acyclically color subcubic graphs with four colors. 

As an additional remark, coloring of subcubic graphs has been an interesting research subject, see for example~ \cite{BS08, DK08, DST08, JMS98, P06, S04}. 

\section{Proof of Theorem \ref{main}}
For every vertex $v$, denote the degree of $v$ by $\deg(v)$, and denote the set of neighbors of $v$ by $N(v)$, and write $N(v) \cup \{v\}$ as $N[v]$.
The {\it order} of a graph is the number of its vertices.

Theorem 2 holds for graphs of order at most $4$ since we can color each vertex by a different color.
We say that $G$ is a {\it minimum counterexample} if it is a subcubic graph without $K_{3,3}$ or $C_5$ as components and there is a family $L$ of lists of size $4$  such that $G$ has no superlinear $L$-coloring, but every subcubic graph with fewer vertices than $G$ and with no component isomorphic to $K_{3,3}$ or $C_5$ is superlinearly $4$-choosable.
Thus every minimum counterexample has at least five vertices.

In the rest of this section, we assume that $G$ is a minimum counterexample and $L$ is a family of lists of size four such that $G$ has no superlinear $L$-coloring.

\begin{lemma}  \label{not a cycle}
$G$ is a connected graph of minimum degree at least two, and $G$ is not a cycle.
\end{lemma}

\begin{pf}
If $G$ is not connected, then $G$ has a component that is not superlinearly $4$-choosable, contradicting the assumption that $G$ is a minimum counterexample.
Suppose that $G$ contains a vertex $v$ of degree one, and let $u$ be the neighbor of $v$ in $G$.
Then $G-v$ is either $C_5$ or a graph that has a superlinear $L$-coloring $f$.
For the latter, we can extend $f$ to $G$ by coloring $v$ with a color different from $f(u)$ and the colors on the neighbors of $u$.
For the former and the case that the $C_5$ does not have a superlinear $L$-coloring, we can first define $f$ on the $C_5$ such that the two neighbors of $u$ receive the same color, but all other vertices in $C_5$ receive the different colors, and then define $f(v)$ to be the color that is different from $f(u)$ and the colors of neighbors of $u$.
So $G$ has minimum degree at least two.

If $G$ is a cycle, then $G$ has at least six vertices, and $G^2$ has maximum degree four but is not $K_5$, where $G^2$ is the graph that is obtained from $G$ by adding edges $uv$ for every pair of vertices $u$ and $v$ with distance two.
So, $G^2$ is $4$-choosable by Theorem \ref{list Brooks}.
Since every proper $L$-coloring of $G^2$ is a superlinear $L$-coloring of $G$, we have that $G$ is superlinearly $4$-choosable.
\end{pf}

\begin{lemma}  \label{delete deg 2 path}
$G$ contains no adjacent veritces of degree two.
\end{lemma}

\begin{pf}
Suppose that $G$ contains adjacent vertices of degree two.
Let $P=v_1v_2...v_k$ be a maximal induced path of order at least two in $G$ such that $\deg(v_i)=2$ for $1 \leq i \leq k$.
Let $u_1$ and $u_k$ be the neighbor of $v_1$ and $v_k$ other than $v_2$ and $v_{k-1}$, respectively.
As $P$ is maximal, $\deg(u_1)=\deg(u_k)=3$.
We claim that $G-P$ is the disjoint union of two $5$-cycles.

Note that $G-P$ has no component isomorphic to $K_{3,3}$ as $G$ does not.
Suppose that $G-P$ contains at most one component isomorphic to $C_5$.
Furthermore, we may assume that if $G-P$ has such a component, then it contains $u_1$.
By induction, there is an $L$-coloring $f$ defined on $G-P$ that is superlinear, except that the two neighbors of $u_1$ may receive the same color.
We define $v_0$ to be $u_1$.
Define $f(v_1)$ to be a color in $L(v_1)$ but different from $f(u_1)$ and the $f$-values of the other neighbors of $u_1$, and then define $f(v_i)$ to be a color in $L(v_i) - \{f(v_{i-2}),f(v_{i-1}), f(u_k)\}$ for $2 \leq i \leq k$.
It is clear that the neighbors of every vertex of degree two have different colors, so no 2-colored cycle passes vertices in $P$, and hence there exists no 2-colored cycle.
In addition, note that $u_k$ has degree three, and $u_k$ is adjacent to two vertices unless $u_k=u_1$, so it is clear that no vertex is adjcent to three vertices that have the same color.
That is, $f$ is a superlinear $L$-coloring of $G$, contradicting the assumption that $G$ is a counterexample.
Therefore, $G-P$ contains two components, and each of them is a $5$-cycle.

As a result, $G$ can be obtained from two disjoint $5$-cycles by adding a path connecting them.
However, there exists a maximal path $P'$ of order four in a $5$-cycle such that every vertex of $P'$ is of degree two, but $G-P'$ does not contain two components, contradicting the claim.
This proves that $G$ contains no adjacent vertices of degree two.
\end{pf}

\bigskip

We say that an induced path $P=v_1v_2...v_k$ in $G$ is {\it special} if $k \geq 3$, and $\deg(v_1)=\deg(v_k)=2$ and $\deg(v_i)=3$ for $2 \leq i \leq k-1$, and the neighbor of $v_1$ (and $v_k$, respectively) outside $P$, denoted by $u_1$ (and $u_k$, respectively) has degree three.
For $2 \leq i \leq k-1$, consider the neighbor of $v_i$ other than $v_{i-1}$ and $v_{i+1}$.
We denote the said neighbor by $u_i$ if it has degree three; otherwise, we denote it by $x_i$, and we denote the neighbor of $x_i$ other than $v_i$ by $u_i$.
For the former case, $x_i$ is undefined.
Note that when $x_i=x_j$ for some $2 \leq i < j \leq k-1$, then $u_i=v_j$ and $u_j=v_i$.

\begin{lemma}  \label{delete path}
Let $P=v_1v_2...v_k$ be a special path in $G$, and define $x_i$ and $u_i$ as in the last paragraph.
Let $Q$ be the subgraph of $G$ induced by $V(P) \cup \{x_i: 2 \leq i \leq k-1$ and $x_i$ is defined$\}$.
If $x_i \neq x_j$ for $2 \leq i <j \leq k-1$, and $u_1$ is not adjacent to three pairwise nonadjacent vertices of $Q$, and no $u_i$ is adjacent to two nonadjacent vertices of $Q$ unless $u_i=u_1$, then there is a component of $G-Q$ isomorphic to $C_5$ containing $u_j$ for some $2 \leq j \leq k$ but not containing $u_1$.
\end{lemma}

\begin{pf}
Note that $G-Q$ has no component isomorphic to $K_{3,3}$ as $G$ does not.
We may assume that $G-Q$ has at most one component isomorphic to $C_5$, for otherwise the lemma holds.
Furthermore, we may assume that if $G-Q$ has such a component, then it contains $u_1$.
By induction, there is an $L$-coloring $f$ defined on $G-Q$ that is superlinear, except that the two neighbors of $u_1$ may receive the same color.
Note that every two vertice in $Q-P$ have degree two, so they are not adjacent, and $\deg(u_i)=3$ for $1 \leq i \leq k$, by Lemma \ref{delete deg 2 path}.

Define $f(v_1)$ to be a color in $L(v_1)$ but different from $f(u_1)$ and the $f$-values of the other neighbors of $u_1$, and then define $f(v_i)$ to be a color in $L(v_i) - \{f(v_{i-1}), f(u_{i-1}), f(u_i)\}$ for $2 \leq i \leq k-2$, where $f(v_{k-2})$ is chosen such that $L(v_{k-1})-\{f(v_{k-2}), f(u_{k-2}), f(u_{k-1}), f(u_k)\}$ is as large as possible.
If $f(v_{k-2})$, $f(u_{k-1})$, $f(u_k)$ are pairwise distinct, then we define $f(v_{k-1})$ to be a color in $L(v_{k-1}) - \{f(v_{k-2}), f(u_{k-1}), \allowbreak f(u_k)\}$; otherwise, we define $f(v_{k-1})$ to be $L(v_{k-1}) - \{f(v_{k-2}), f(u_{k-2}), \allowbreak f(u_{k-1}), \allowbreak f(u_k)\}$.
And finally define $f(v_k)$ to be a color in $L(v_k) - \{f(v_{k-2}), \allowbreak f(v_{k-1}), \allowbreak f(u_k)\}$, and define $f(x_i)$ to be a color in $L(x_i) - \{f(v_{i-1}), f(v_i), \allowbreak f(u_i)\}$ for those $i$ such that $x_i$ exists.

It is clear that the neighbors of every vertex of degree two have different colors, so no two-colored cycle passes vertex of degree two.
Suppose there is a two-colored cycle, then it must pass some vertex $v_i$ for $2 \leq i \leq k-1$.
Let $t$ be the minimum number such that $v_t$ is in a 2-colored cycle, so $u_t$ is adjacent to $v_t$, and $u_t$, $v_{t+1}$ are also in the 2-colored cycle, and $t \leq k-2$.
However, it is impossible since $f(v_{t+1}) \neq f(u_t)$ when $t+1 \leq k-2$.
Also, if $t=k-2$, then $u_{k-1}$ is adjacent to $v_{k-1}$, and $u_{k-1}$ is also in the 2-colored cycle.
However, it is still impossible since either $f(v_{k-1}) \neq f(u_{k-2})$ or $f(v_{k-2}) \neq f(u_{k-1})$.
Hence there are no 2-colored cycles.

In addition, it is clear that no vertex other than $v_{k-2}$ is adjcent to three vertices that have the same color.
Suppose that $v_{k-2}$ is adjacent to three vertices that have the same color, then $f(v_{k-3})=f(u_{k-2})$ and the three colors $f(v_{k-2})$, $f(u_{k-1})$, $f(u_k)$ are pairwise distinct.
This implies that $L(v_{k-2}) - \{f(v_{k-3}), f(u_{k-3}), f(u_{k-2})\}$ has size two, and $L(v_{k-1}) = \{f(u_{k-2}), f(v_{k-2}),f(u_{k-1}),f(u_k)\}$.
Hence, if we choose the other color in $L(v_i) - \{f(v_{i-1}), f(u_{i-1}), f(u_i)\}$ when we color $v_{k-2}$, $L(v_{k-1})-\{f(v_{k-2}), f(u_{k-2}), f(u_{k-1}), f(u_k)\}$ is larger, a contradiction.
Consequently, $f$ is a superlinear $L$-coloring of $G$.
\end{pf}

\begin{lemma} \label{eyeglass}
Let $Q=v_1v_2v_3v_4$ be an induced path in $G$ with $\deg(v_1)=\deg(v_4)=2$, and let $u_1$ and $u_4$ be the neighbor of $v_1$ and $v_4$ outside $Q$, respectively.
If $u_1$ is adjacent to $v_3$, then $u_4$ is not adjacent to $v_2$.
\end{lemma}

\begin{pf}
Suppose that $u_1$ is adjacent to $v_3$, and $u_2$ is adjacent to $v_4$.
Since $v_1$ and $v_4$ have degree two, $u_1$ and $u_4$ are of degree three by Lemma \ref{delete deg 2 path}, so $G-Q$ does not contain a component isomorphic to $C_5$ or $K_{3,3}$.
Therefore, there exists a superlinear $L$-coloring of $G-Q$.
For $i=2,3$, define $f(v_i)$ to be a color in $L(v_i)-\{f(u_1),f(u_4)\}$ such that $f(v_2) \neq f(v_3)$.
And then for $i=1,4$, define $f(v_i)$ to be a color in $L(v_i)-\{f(u_i),f(v_2),f(v_3)\}$.
It is clear that $f$ is a superlinear $L$-coloring of $G$, a contradiction.
\end{pf}

\begin{lemma} \label{C4 2 deg 2}
Let $P=v_1v_2v_3v_4v_5$ be an induced path in $G$ with $\deg(v_1)=\deg(v_5)=2$, and let $u_i$ be the neighbor of $v_i$ outside $P$ for $i=1,3,5$.
Let $w_3$ be a common neighbor of $v_2$ and $v_4$ other than $v_3$, and let $x$ be the neighbor of $w_3$ other than $v_2$ and $v_4$.
Let $Q$ be the subgraph of $G$ induced by $V(P) \cup \{w_3\}$.
If $\deg(x)=\deg(u_3)=3$, and $w_3 \not \in \{u_1,u_5\}$, then there exists a component of $G-Q$ isomorphic to $C_5$ not containing $u_1$ or $u_5$.
\end{lemma}

\begin{pf}
Suppose that either $G-Q$ contains no component isomorphic to $C_5$, or each component of $G-Q$ isomorphic to $C_5$ contains $u_1$ or $u_5$.
By induction, there is an $L$-coloring $f$ defined on $G-Q$ that is superlinear, except that the colored neighbors of $u_1$ may receive the same color, and the colored neighbors of $u_5$ may receive the same color.

We shall consider two cases.
The first case is that $(L(v_2)-\{f(u_1),\allowbreak f(u_3)\}) \cap (L(v_4)-\{f(u_5),f(x)\}) \neq \emptyset$.
In this case, define $f(v_2)=f(v_4)$ to be a color in $(L(v_2)-\{f(u_1),f(u_3)\}) \cap (L(v_4)-\{f(u_5),f(x)\})$.
And then define $f(v_1)$ to be a color in $L(v_1)-\{f(u_1),f(v_2)\}$ but different from the $f$-value of some neighbor of $u_1$, and define $f(v_5)$ to be a color in $L(v_5)-\{f(v_4),f(u_5)\}$ but different from the $f$-value of some neighbor of $u_5$.
Finally, define $f(v_3)$ to be a color in $L(v_3)-\{f(v_2),f(u_3)\}$ but different from the $f$-value of some neighbor of $u_3$, and define $f(w_3)$ to be a color $L(w_3)-\{f(v_2), f(v_3),f(x)\}$.

The second case is that $(L(v_2)-\{f(u_1),f(u_3)\}) \cap (L(v_4)-\{f(u_5),\allowbreak f(x)\}) = \emptyset$.
Define $f(v_3)$ to be a color in $L(v_3)-\{f(u_3)\}$ but different from the $f$-value of some neighbor of $u_3$.
Since $(L(v_2)-\{f(u_1),f(u_3)\}) \cap (L(v_4)-\{f(u_5),f(x)\}) = \emptyset$, without loss of generality, we may assume that $f(v_3)$ is not in $f(v_4)-\{f(u_5),f(x)\}$.
Then, define $f(v_2)$ to be a color in $L(v_2)-\{f(u_1),f(u_3),f(v_3)\}$, and define $f(v_1)$ to be a color in $L(v_1)-\{f(u_1), f(v_2)\}$ but different from an $f$-value of some neighbor of $u_1$, and define $f(w_3)$ to be a color in $L(w_3)-\{f(v_2),f(v_3),f(x)\}$.
Finally, define $f(v_4)$ to be a color in $L(v_4)-\{f(w_3), f(u_5),f(x)\}$, and define $f(v_5)$ to be a color in $L(v_5)-\{f(v_4),f(u_5)\}$ but different from an $f$-value of some neighbor of $u_5$.
It is clear that $f$ is a superlinear $L$-coloring in the both cases, a contradiction, so $G-Q$ contains a component isomorphic to $C_5$ not containing $u_1$ or $u_5$.
\end{pf}

\begin{lemma}  \label{distance}
The distance between any two vertices of $G$ of degree two is at least five.
\end{lemma}

\begin{pf}
Let $a$ and $b$ be two different vertices of degree two such that the distance between $a$ and $b$ is as small as possible.
Let $P=v_1v_2...v_{t+1}$ be a shortest path from $a$ to $b$, where $a=v_1$, $b=v_{t+1}$ and $t$ is the length of $P$.
Since there are no adjacent vertices of degree two by Lemma \ref{delete deg 2 path}, $P$ is a special path, and we define vertices $u_i,x_i$ and graph $Q$ as in Lemma \ref{delete path}.

Suppose that $t=2$.
Since $G$ contains at least five vertices, $G=K_{2,3}$ if $u_1$ or $u_k$ is adjacent to three vertices of $Q$.
However, it is impossible since $K_{2,3}$ is superlinearly $4$-choosable: assigning two different colors $c_1,c_2$ to the two vertices of degree three in $K_{2,3}$, and assigning colors different from $c_1,c_2$ to the three vertices of degree two such that these three vertices do not receive the same color.
On the other hand, if $u_i$ is adjacent to two nonadjacent vertices of $Q$ for some $i \geq 2$, then $u_1 \neq u_2=u_3$, so we may assume that $u_i$ is not adjacent to two nonadjacent vertices of $Q$ unless $u_i=u_1$, by swapping $v_1$ and $v_3$. 
By Lemma \ref{delete path}, $G-Q$ contains a 5-cycle as a component which does not contain $u_1$, so there are at most two edges with one end in $Q$ and one end in the component.
This implies that two vertices of degree two in the 5-cycle are adjacent, a contradiction.

Hence, $t \geq 3$, and $x_i$ does not exist for $2 \leq i \leq t$ as $v_1,v_{t+1}$ is the closest pair of vertices of degree two.
Similarly, $u_i=u_j$ implies that $\lvert j-i \rvert \leq 2$, so no $u_i$ is adjacent to three pairwise nonadjacent vertices of $Q$.

If $t=3$, then by Lemma \ref{eyeglass} and by symmetry, we may assume that $u_i$ is adjacent to two nonadjacent vertices of $Q$ only when $u_i=u_1$.
Therefore, by Lemma \ref{delete path}, there exists a component of $G-Q$ isomorphic to $C_5$ but not containing $u_1$.
So there are at most three edges with one end in $Q$ and one end in the $5$-cycle, and hence there exist two vertices of degree two with distance $2$, a contradiction.
Similarly, when $t=4$, $u_2 \neq u_4$ or by Lemma \ref{C4 2 deg 2}, $G-Q$ contains a $5$-cycle $C$ as a component which does not contain $u_1$, so $C$ contains a vertex $v$ of degree two in $G$.
If $v_5$ is adjacent to some vertex of $C$, then the distance between $v$ and $v_5$ is smaller than $t$, a contradiction.
If $v_5$ is not adjacent to any vertex of $C$, then there are only three edges with one end in $Q$ and one end in $C$, so $C$ contains two vertices of degree two between distance less than $t$, a contradiction.
This proves that $t \geq 5$.
\end{pf}

\begin{lemma}  \label{triangle-free}
$G$ contains no triangle as a subgraph.
\end{lemma}

\begin{pf}
Let $C=v_1v_2v_3v_1$ be a triangle in $G$.
By Lemma \ref{distance}, there is at most one $v_i$ of degree two, so we may assume that $3 = \deg(v_1)=\deg(v_2) \geq \deg(v_3)$.
For $i=1,2,3$, let $x_i$ be the neighbor of $v_i$ other than vertices in $C$ such that $\deg(x_i)=2$ if such vertex exists.
Let $Q = \{v_i,x_i: 1 \leq i \leq 3\}$ (we only consider those $x_i$ which are defined).
Let $u_i$ be the neighbor of $x_i$ other than vertices in $Q$ if $x_i$ is defined, and let $u_i$ be the neighbor of $v_i$ other than vertices in $Q$ if $x_i$ is not defined and $\deg(v_i)=3$.
Note that $G-Q$ contains no $K_{3,3}$ as a component.
Since the distance of any two vertices of degree two is at least five, and the number of edges which have one end in $Q$ and the other end in $G-Q$ is at most three, $G-Q$ contains no 5-cycle as a component.
Hence, we can apply induction to $G-Q$ to obtain a superlinear $L$-coloring $f$ of $G-Q$.
Now, we shall extend $f$ to a superlinear $L$-coloring of $G$.

Note that we ignore $f(u_3)$ in the following sentence if $u_3$ is not defined.
Define $f(v_1)$ to be a color in $L(v_1) - \{f(u_1), f(u_2), f(u_3)\}$, $f(v_2)$ to be a color in $L(v_2) - \{f(v_1), f(u_2), f(u_3)\}$, and $f(v_3)$ to be a color in $L(v_3) - \{f(v_1), f(v_2), f(u_3)\}$, and then define $f(x_i)$ to be a color in $L(x_i) - \{f(v_i), f(u_i)\}$ for those $i$ such that $x_i$ are defined.
It is clear that the neighbors of any vertex of degree two receive different colors, so no 2-colored cycles pass through a vertex of degree two.
If there is a 2-colored cycle, then it must contain the path $u_i v_i v_j u_j$ for some $i<j$, and $u_i$ is adjacent to $v_i$, and $u_j$ is adjacent to $v_j$, but that is impossible since $f(v_i) \neq f(u_j)$.
Therefore, there are no 2-colored cycles.
In addition, as the neighbors of $u_i$ (other than $v_i$) have different colors, it is clear that no vertex is adjacent to three vertices that have the same color, so $f$ is a superlinear $L$-coloring of $G$.
\end{pf}

\begin{lemma} \label{C4 contains deg 2}
No $4$-cycle in $G$ contains a vertex of degree two.
\end{lemma}

\begin{pf}
Suppose that $G$ contains a $4$-cycle $Q=v_1v_2v_3v_4v_1$ with $\deg(v_1)=2$.
Let $u_i$ be the neighbor of $v_i$ outside $C$ for $2 \leq i \leq 4$.
By Lemma \ref{distance}, $\deg(v_i)=\deg(u_i)=3$ for $2 \leq i \leq 4$, and $G-Q$ does not contain $K_{3,3}$ or $C_5$ as a component.
Furthermore, if some vertex outside $Q$ is adjacent to two vertices of $Q$, then $u_2=u_4$ is the only such vertex by Lemma \ref{triangle-free}.
By induction, $G-Q$ has a superlinear $L$-coloring $f$.
Define $f(v_3)$ to be a color in $L(v_3)-\{f(u_2),f(u_3),f(u_4)\}$, $f(v_2)$ to be a color in $L(v_2)-\{f(u_2),f(v_3)\}$ but different from the $f$-value of a colored neighbor of $u_2$, and define $f(v_4)$ to be a color in $L(v_4)-\{f(v_2),f(v_3),f(u_4)\}$, $f(v_1)$ to be a color in $L(v_1)-\{f(v_2),f(v_3),f(v_4)\}$.
Clearly, $f$ is a superlinear $L$-coloring of $G$, a contradiction.
\end{pf}

\begin{lemma} \label{C4 adjacent to deg 2}
No vertex of degree two in $G$ is adjacent to a $4$-cycle.
\end{lemma}

\begin{pf}
Suppose that there exist a $4$-cycle $C=v_1v_2v_3v_4v_1$ in $G$ and a vertex $v$ of degree two adjacent to $v_1$.
Let $u_1$ be the neighbor of $v$ other than $v_1$, and $u_i$ the neighbor of $v_i$ outside $C$ for $2 \leq i \leq 4$.
Let $Q$ be the subgraph of $G$ induced by $\{v,v_i:1 \leq i \leq 4\}$.
Note that $v \neq u_3$ by Lemma \ref{C4 contains deg 2}.
By Lemma \ref{distance}, $\deg(u_i)=3$ for $1 \leq i \leq 4$, and $G-Q$ contains no component isomorphic to $K_{3,3}$ or $C_5$.
Furthermore, Lemmas \ref{triangle-free} and \ref{C4 contains deg 2} ensure that $\{u_1,u_3\} \cap \{u_2,u_4\}=\emptyset$.
By induction, $G-Q$ has a superlinear $L$-coloring $f$.
Define $f(v_3)$ to be a color in $L(v_3)-\{f(u_2),f(u_3),f(u_4)\}$, $f(v_2)$ to be a color in $L(v_2)-\{f(u_2),f(v_3)\}$ but different from the $f$-value of a neighbor of $u_2$, $f(v_4)$ to be a color in $L(v_4)-\{f(v_2),f(v_3),f(u_4)\}$, $f(v_1)$ to be a color in $L(v_1)-\{f(u_1),f(v_2),f(v_4)\}$, and define $f(v)$ to be a color in $L(v)-\{f(u_1),f(v_1)\}$ but different from the $f$-vaule of a neighbor of $u_1$.
Clearly, $f$ is a superlinear $L$-coloring of $G$, a contradiction.
\end{pf}

\begin{lemma} \label{C5 contains deg 2}
No $5$-cycle in $G$ contains a vertex of degree two.
\end{lemma}

\begin{pf}
Suppose that there exist a $5$-cycle $Q=v_1v_2v_3v_4v_5v_1$ in $G$ with $\deg(v_1)=2$.
Let $u_i$ be the neighbor of $v_i$ outside $Q$ for $2 \leq i \leq 5$.
By Lemma \ref{distance}, $\deg(v_i)=\deg(u_i)=3$ for $2 \leq i \leq 5$, and $G-Q$ contains no component isomorphic to $K_{3,3}$ or $C_5$.
Furthermore, Lemmas \ref{triangle-free}, \ref{C4 contains deg 2} and \ref{C4 adjacent to deg 2} ensure that $u_i$ are pairwise distinct for $2 \leq i \leq 5$.
By induction, $G-Q$ has a superlinear $L$-coloring $f$.
Define $f(v_3)$ to be a color in $L(v_3)-\{f(u_2),f(u_3),f(u_4)\}$, $f(v_4)$ to be a color in $L(v_4)-\{f(v_3),f(u_4),f(u_5)\}$, $f(v_2)$ to be a color in $L(v_2)-\{f(u_2),f(v_3),f(v_4)\}$, $f(v_5)$ to be a color in $L(v_5)-\{f(v_2),f(v_4),f(u_5)\}$, and define $f(v_1)$ to be a color in $L(v_1)-\{f(v_2), \allowbreak f(v_5)\}$.
Clearly, $f$ is a superlinear $L$-coloring of $G$, a contradiction.
\end{pf}

\begin{lemma}\label{cubic}
$G$ is cubic. 
\end{lemma}

\begin{pf}
Let $P=v_1v_2v_3v_4$ be a path in $G$, where the degree of $v_4$ is two in $G$.
Since the distance between any two vertices of degree two is at least five, $v_i$ and its neighbor that is not in $P$ are of degree three for $1 \leq i \leq 4$.
Furthermore, $P$ is an induced path, and the neighbors of each $v_i$ are pairewise distinct by Lemmas \ref{triangle-free}, \ref{C4 contains deg 2}, \ref{C4 adjacent to deg 2} and \ref{C5 contains deg 2}.
For each $1 \leq i \leq 4$, let $u_i$ be a neighbor of $v_i$ other than vertices in $P$, and let $w$ be the neighbor of $v_1$ other than $u_1$ and $v_2$.
Note that the number of edges which have one end in $P$ and the other end in $G-P$ is at most five, and the distance of any two vertices of degree two is at least five, so $G-P$ does not contain $C_5$ as a component when there are at most four edges between $P$ and $G-P$.
On the other hand, if there are exactly five edges between $P$ and $G-P$ such that $G-P$ contains a 5-cycle as a component, then $G$ contains exactly $9$ vertices, and $v_1$ together with three vertices of the $5$-cycle induce a $4$-cycle by Lemma \ref{triangle-free}, so there is a path $P'$ having $v_4$ as an end such that $P'$ is disjoint with the $4$-cycle, and hence $G-P'$ has no component isomorphic to $C_5$, and we replace $P$ by $P'$.
Therefore, we may assume that $G-P$ does not contain $C_5$ as a component.
In addition, $G-P$ does not contain $K_{3,3}$ as a component.
Hence, we can apply induction to $G-P$ to obtain a superlinear $L$-coloring $f$ of $G-P$.
And we shall extend $f$ to a superlinear $L$-coloring of $G$.

Suppose that $f(u_1) \neq f(w)$.
Define $f(v_1)$ to be a color in $L(v_1) - \{f(u_1), f(w), f(u_2)\}$, $f(v_2)$ to be a color in $L(v_2) - \{f(v_1), f(u_2), f(u_3)\}$, $f(v_3)$ to be a color in $L(v_3) - \{f(v_2), f(u_3), f(u_4)\}$, and define $f(v_4)$ to be a color in $L(v_4) - \{f(u_3), f(v_3), f(u_4)\}$.
Note that neighbors of vertices of degree two receive different colors, so no 2-colored cycle passes $v_4$.
If there is a 2-colored cycle, then it must pass through $v_i, v_j, u_j$ for some $i<j<4$ and $v_k$ for all $k$ with $i \leq k \leq j$ since $f(u_1) \neq f(w)$, but it is a contradiction since $f(v_{j-1}) \neq f(u_j)$.
Therefore, there are no 2-colored cycles.
And it is clear that no vertex is adjacent to three vertices that receive the same color, so $f$ is a superlinear $L$-coloring and $G$ cannot be a counterexample, a contradiction.
So, let $f(u_1)=f(w)$.

Define $f(v_1)$ to be a color in $L(v_1) - \{f(z): z \in N[w] - \{v_1\}\}$.
If $f(v_1) = f(u_2)$, then define $f(v_3)$ to be a color in $L(v_3) - \{f(u_2), f(u_3), \allowbreak f(u_4)\}$, $f(v_2)$ to be a color in $L(v_2) - \{f(u_1), f(u_2), f(v_3)\}$, and $f(v_4)$ to be a color in $L(v_4) - \{f(u_3), f(v_3), f(u_4)\}$.
If $f(v_1) \neq f(u_2)$, then define $f(v_2)$ to be a color in $L(v_2) - \{f(u_1), f(v_1), f(u_2)\}$, and further define $f(v_3)$ to be a color in $L(v_3) - \{f(u_2), f(v_2), f(u_4)\}$ ($L(v_3) - \{f(v_2),f(u_3),f(u_4)\}$, respectively) when $f(v_2)=f(u_3)$ ($f(v_2) \neq f(u_3)$, respectively), and define $f(v_4)$ to be a color in $L(v_4) - \{f(u_3), f(v_3), \allowbreak f(u_4)\}$.
It is clear that neighbors of a vertex of degree two receive different colors.
Also, if a 2-colored cycle exists, then it must pass through $u_i,v_i,v_j, u_j$ for some $i<j<4$ and $v_k$ for $i \leq k \leq j$, but it is impossible since $f(v_2) \neq f(u_1)$ and either $f(v_3) \neq f(u_2)$ or $f(u_3) \neq f(v_2)$.
And no vertex is adjacent to three vertices that receive the same color.
Hence, $f$ is a superlinear $L$-coloring of $G$, contradicting that $G$ is a counterexample.
This completes the proof.~
\end{pf}

\begin{lemma}  \label{K_{2,3}-free}
$G$ does not contain $K_{2,3}$ as a subgraph.
\end{lemma}

\begin{pf}
Let $H$ be a subgraph of $G$ isomorphic to $K_{2,3}$.
In fact, $H$ is an induced subgraph by Lemma \ref{triangle-free}.
Let $V(H) = \{v_1,v_2, u_1,u_2,u_3\}$ and $E(H) = \{v_iu_j: i=1,2, j=1,2,3\}$.
Note that $v_1$ and $v_2$ are all the common neighbors of $u_1,u_2,u_3$, since $G$ does not contain $K_{3,3}$ as a component.

Suppose that $u_1$ and $u_2$ have a common neighbor $v_3$ other than $v_1$ and $v_2$.
Note that $G-(H \cup \{v_3\})$ does not contain $K_{3,3}$ or $C_5$ as a component by Lemma \ref{distance}.
So we can apply induction to $G-(H \cup \{v_3\})$ to obtain a superlinear $L$-coloring $f$ of $G-(H \cup \{v_3\})$, and we shall extend $f$ to $G$.
Let $x$ be the neighbor of $u_3$ outside $H$, and $y$ the neighbor of $v_3$ outside $H$.
Define $f(u_3)$ to be a color in $L(u_3) - \{f(v): v \in N[x] - \{u_3\}\}$, $f(v_3)$ to be a color in $L(v_3) - \{f(v): v \in N[y] - \{v_3\}\}$, $f(u_2)$ to be a color in $L(u_2) - \{f(u_3), f(v_3), f(y)\}$, $f(v_2)$ to be a color in $L(v_2) - \{f(u_2), f(u_3), f(x)\}$, $f(v_1)$ to be a color in $L(v_1) - \{f(u_2), f(v_2), f(u_3)\}$, and $f(u_1)$ to be a color in $L(u_1) - \{f(v_1), f(v_2), f(v_3)\}$ ($L(u_1) - \{f(v_1), f(v_2), f(v_3), f(u_2)\}$, respectively) when $f(v_1), f(v_2), f(v_3)$ are pairwise distinct (when $f(v_1)$,$f(v_2)$,$f(v_3)$ are not pairwise distinct, respectively).
By the choices of colors,  we see that no vertex is adjacent to three vertices of the same color and there are no two-colored cycles, so $G$ is not a counterexample. 
Hence, we may assume that $v_1$ and $v_2$ are the only common neighbors of any two of $u_1,u_2,u_3$ by symmetry.
In other words, it is impossible to add an edge to $G$ to make $G$ contain $K_{3,3}$ as a subgraph.

Let $a,b,c$ be the three vertices in $G-H$ adjacent to $u_1, u_2, u_3$ in $G$, respectively. 
Note that $G-H$ does not contain 5-cycles and $K_{3,3}$ as components, so we can apply induction to $G-H$ and obtain a superlinear $L$-coloring $f$ of $G-H$.   
Now we extend this coloring to $G$. 

We define $f(u_1), f(u_2)$ so that $f(u_1)$ is different from $f(a)$ and the colors of the two other neighbors of $a$, and $f(u_2)$ is different from $f(u_1), f(b)$.
Now define $f(v_1)$ so that $f(v_1)\in L(v_1)-\{f(u_1), f(u_2), \allowbreak f(c)\}$ and define $f(v_2)$ so that $f(v_2)\in L(v_2)-\{f(v_1), f(u_1), f(u_2)\}$.  Finally we define $f(u_3)$ so that $f(u_3)\in L(u_3)-\{f(c), f(v_1), f(v_2)\}$ and if possible, $f(u_3)\not=f(u_2)$.   

By induction, the two neighbors of $c$ other than $u_3$ get different colors, $c$ is not adjacent to three vertices of the same color.  Since $f(v_1)\not=f(c)$, $u_3$ is not adjacent to three vertices of the same color.  Since $u_1, u_2, v_1, v_2$ have distinct colors,  none of them is adjacent to three vertices of the same color.  By the choices of $f(u_1)$ and $f(u_2)$, none of $a$ and $b$ is adjacent to three vertices of the same color.  

Since $u_1, u_2, v_1, v_2$ have distinct colors, no 2-colored cycle is in $H$.  Because of the choice of $f(u_1)$,   no 2-colored cycle contains $a$ and $u_1$.   By the choice of $f(u_3)$, either $f(u_3)\not=f(u_2)$, or $f(c)=f(u_2)\not=f(v_2)$ (and $f(v_1)\not=f(c)$),   so no $2$-colored cycle contains $b$ and $c$. 

So the extension gives a superlinear $L$-coloring of $G$. 
\end{pf}

\bigskip

Now we prove a lemma about list-coloring.

\begin{lemma}  \label{list-color cycle}
Let $C$ be a cycle and $L'=\{L'(v): v \in V(G)\}$ a family of lists such that $\lvert L'(v) \rvert \geq 2$ for every vertex $v$ of $C$.
If either $\lvert L'(w) \rvert \geq 3$ or $L'(u) \neq L'(v)$ for some vertices $u,v,w$ of $C$, then $C$ has a linear $L'$-coloring.
\end{lemma}

\begin{pf}
Let $C=v_1v_2...v_kv_1$ and either $\lvert L'(v_k) \rvert \geq 3$ or $L'(v_1) \neq L'(v_k)$.
Define an $L'$-coloring $f$ on $C$ by letting $f(v_1)$ be a color in $L'(v_1)-L'(v_k)$ if possible, and define $f(v_i)$ to be a color in $L'(v_i) - \{f(v_{i-1})\}$ for $2 \leq i \leq k-1$, and define $f(v_k)$ to be a color in $L'(v_k) - \{f(v_1), f(v_{k-1})\}$.
Notice that $f$ is a proper $L'$-coloring since $L'(v_k) - \{f(v_1), f(v_{k-1})\}$ is not empty by the choice of $f(v_1)$.
If $C$ is 2-colored, then $f(v_1)=f(v_{k-1})$, and we can redefine $f(v_k)$ to be a color in $L'(v_k) - \{f(v_1)=f(v_{k-1}), f(v_{k-2})\}$ to make $f$ be a linear $L'$-coloring in this case.
\end{pf}

\begin{lemma}  \label{good cycle}
If $C$ is an induced cycle in $G$ such that no vertex is adjacent to at least two vertices in $C$, then $G-C$ contains $C_5$ as a component.
\end{lemma}

\begin{pf}
Let $C=v_0v_1...v_{k-1}v_0$ and $u_i$ be the neighbor of $v_i$ other than vertices in $C$ for $0 \leq i \leq k-1$.
Suppose that $G-C$ does not contain $C_5$ as a component, so we can apply induction to $G-C$ to obtain a superlinear $L$-coloring $f$ of $G-C$.
Define $L^+(v_i) = L(v_i) - \{f(u_i), f(u_{i+1})\}$ and $L^-(v_i) = L(v_i) - \{f(u_i), f(u_{i-1})\}$ for $0 \leq i \leq k-1$, where the indices are computed modulo $k$.
Note that if $C$ has a proper $L^+$-coloring or a proper $L^-$-coloring, then this coloring together with $f$ gives a proper $L$-coloring of $G$ such that no vertex is adjacent to three vertices of the same color, and the only possible 2-colored cycle is $C$ since no 2-colored path contains some vertices in $C$ as internal vertices.
As a result, $f$ can be extended to a superlinear $L$-coloring if $C$ has a linear $L^+$-coloring or a linear $L^-$-coloring.

Suppose that $C$ does not have a linear $L^+$-coloring nor a linear $L^-$-coloring.
By Lemma \ref{list-color cycle}, $L^+(v_i) = L^+(v_{i+1})$, $L^-(v_i) = L^-(v_{i+1})$, $\lvert L^+(v_i) \rvert = \lvert L^-(v_i) \rvert =2$, for every $0 \leq i \leq k-1$.
So $f(u_{i-1}) \neq f(u_i) \neq f(u_{i+1})$ for all $0 \leq i \leq k-1$, for otherwise $L^+(v_i)$ or $L^-(v_i)$ has size at least three. 
Furthermore, $L(v_i) = L^+(v_i) \cup \{f(u_i), f(u_{i+1})\}$ and $L(v_{i+1}) = L^+(v_{i+1}) \cup \{f(u_{i+1}), f(u_{i+2})\}$, so $L^-(v_i) = L^+(v_i) \cup \{f(u_{i+1})\} - \{f(u_{i-1})\}$ and $L^-(v_{i+1}) = L^+(v_{i+1}) \cup \{f(u_{i+2})\} - \{f(u_i)\}$ for every $0 \leq i \leq k-1$.
Since $L^-(v_i) = L^-(v_{i+1})$ and $f(u_{i-1}) \neq f(u_i) \neq f(u_{i+1})$, $f(u_{i+1}) = f(u_{i-1})$ for every $0 \leq i \leq k-1$.
If $C$ is odd, then $f(u_1)=f(u_2)$, a contradiction.
So $C$ is an even cycle.
In addition, every even cycle is 2-choosable, so $C$ must be 2-colored by every proper $L^+$(and $L^-$)-coloring of $C$.
Without loss of generality, we may assume that $C$ is colored by color $1$ and $2$, and $L^+(v_i)=L^-(v_i)=\{1,2\}$.
So $L(v_i) = \{1,2,f(u_0), f(u_1)\}$ for every $0 \leq i \leq k-1$.

For any $0 \leq i \leq k-1$, we can redefine $f(v_i)$ by a color in $L(v_i)-\{f(u_i), f(v_{i-1})=f(v_{i+1})\}$ and different from the current color, and this alteration will make $f$ be a superlinear $L$-coloring unless the subpath $x_{i-1}u_{i-1}v_{i-1}v_iv_{i+1}u_{i+1}x_{i+1}$ is contained in a 2-colored cycle, where $x_{i-1}$ and $x_{i+1}$ are neighbors of $u_{i-1}$ and $u_{i+1}$ other than vertices in $C$, respectively.
Hence, one neighbor of $u_{i-1}$ and one neighbor of $u_{i+1}$ other than vertices in $C$ has the same $f$-value as $v_{i-1}$, for every $0 \leq i \leq k-1$.
Similarly, we can swap color $1$ and $2$ on $C$ before we change $f(v_i)$, so one neighbors of $u_{i-1}$ has $f$-value $1$ and one neighbor of $u_{i-1}$ has $f$-value $2$ for every $0 \leq i \leq k-1$.
Consequently, let $g$ be the $L$-coloring of $G$ such that $g(v)=f(v)$ for $v \in G-C$, and $g(v_0)=f(u_1)$, $g(v_1)=f(u_0)$, $g(v_{2j})=1$ and $g(v_{2j+1})=2$ for $1 \leq j \leq k/2-1$.
It is easy to see that $g$ is a superlinear $L$-coloring of $G$.
\end{pf}

\bigskip

We are now ready to finish the proof of the theorem.   

\bigskip

\noindent {\bf Proof of Theorem \ref{main}:}
Let $G$ be a minimum counterexample and $C=v_1v_2...v_k$ be a shortest cycle in $G$, so the length of $C$ is at least four by Lemma \ref{triangle-free}.
Suppose that there is a vertex $v$ adjacent to at least two vertices in $C$.
If the shortest cycles have length five or more, then we have a shorter cycle by replacing a segment of $C$ by $v$.
So $C$ has length four, but then $C \cup \{v\}$ induces a $K_{2,3}$, contradicting Lemma \ref{K_{2,3}-free}.
Therefore, there is no vertex adjacent to at least two vertices in $C$, so $G-C$ has a 5-cycle as a component by Lemma \ref{good cycle}.
On the other hand, the shortness of $C$ implies that $C$ has length at most five, so $G$ can be partitioned into two 5-cycles and $G$ has girth five since $G$ is cubic by Lemma \ref{cubic}.
In other words, $G$ is the Petersen graph.

Let $L$ be a family of lists of size four.
Note that if $G-C$ has a superlinear $L$-coloring $f$, then we can extend $f$ to a superlinear $L$-coloring of $G$ by Lemma \ref{good cycle}.
On the other hand, $G-C$ does not have a superlinear $L$-coloring only if all lists $L(v)$ of vertices $v$ in $G-C$ are the same.
This implies that for every 5-cycle $C'$ in $G$, every vertex in $C'$ has the same list of colors.
Hence, every vertex in the Petersen graph has the same list of colors, say $\{1,2,3,4\}$.

We denote the vertex-set of the Petersen graph by $\{u_i,v_i: 1 \leq i \leq 5\}$ such that $u_i$ is adjacent to $v_i$ for $1 \leq i \leq 5$, and $u_1u_2...u_5u_1$ and $v_1v_3v_5v_2v_4v_1$ are $5$-cycles.
Define $f(u_1)=f(u_3)=f(v_5)=1$, $f(u_2)=f(u_4)=2$, $f(v_1)=f(v_2)=f(u_5)=3$, and $f(v_3)=f(v_4)=4$.
Clearly, $f$ is a superlinear $L$-coloring of the Petersen graph, so $G$ is not a counterexample, a contradiction.
$\Box$

\section{Linear-time algorithm}
In this section, we shall give a linear-time algorithm to find a superlinear $L$-coloring whenever a family $L$ of lists is given.

For every vertex $v$ of a graph $G$ and for every positive integer $k$, denote by $N_k(v)$ the set $\{u \in V(G): 0< d(u,v) \leq k\}$, where $d(u,v)$ is the distance between $u$ and $v$, and define $N_k[v] = N_k(v) \cup \{v\}$.

First, we introduce a subroutine to deal with the noncubic case.

\bigskip

\noindent{\bf Subroutine 1}

\noindent{\bf Input:} $(G, L, T)$, where $G$ is a subcubic graph whose every component has minimum degree at most two, $L=\{L(v): v \in V(G)\}$ is a family of lists of size $4$, and $T$ is a list that consists of all vertices of degree at most two in $G$.
Every component of $G$ has minimum degree at most two and does not contain $K_{3,3}$ or $C_5$ as a component.

\noindent{\bf Output:} A superlinear $L$-coloring of $G$.

\noindent{\bf Running time:} $O(\lvert V(G) \rvert)$.

\noindent{\bf Description:}
Pick a vertex $v$ from $T$.
Let $H$ be the component of $G$ containing $v$.
Note that $v$ is of degree at most two in $G$.
If the order of $H$ is at most 10, then output the superlinear $L$-coloring of $G$ by combining a superlinear $L$-coloring of $H$ found by brute force and the superlinear $L$-coloring obtained by executing the subroutine with input $(G-H, L, T-V(H))$.

If $\deg(v) \leq 1$, then put the neighbor of $v$ in $G$ to $T$.
Output the superlinear $L$-coloring of $G$ obtained by extending the superlinear of $G-v$ obtained by executing the subroutine with input $(G-v, L, T-\{v\})$.
Note that $G-\{v\}$ does not contain $C_5$ as a component when $\deg(v)=1$, since $H-\{v\}$ has order at least 9 and $G$ has no component isomorphic to $C_5$.
So we may assume that every vertex in $N_{12}[v]$ has degree at least two.

If there are two adjacent degree two vertices $u,w$ in $N_{12}[v]$, then let $H'$ be the maximal connected subgraph of $H$ containing $u,w$ induced by degree at most two vertices in $H$.
Note that $H'$ has maximum degree at most two, so $H'$ is a cycle or a path.
If $H'$ is a cycle, then $H=H'$, and return the superlinear $L$-coloring obtained by combining a superlinear $L$-coloring of $H'$ and the superlinear $L$-coloring obtained by applying the subroutine with input $(G-H',L,T-V(H'))$.
If $H'$ is a path, then let $G'$ be the subgraph obtained from $G$ by deleting $H'$ and all components of $G-H'$ isomorphic to $C_5$, and let $T'$ be the minimum superset of $T \cap V(G')$ containing $N(H') \cap V(G')$.
Apply the subroutine with input $(G',L,T')$ to obtain a superlinear $L$-coloring of $G'$, and then extend the coloring to a superlinear $L$-coloring of $G$ as in the proof of Lemma \ref{delete deg 2 path}.
Hence, we assume that no two adjacent degree two vertices $u,w$ in $N_{12}[v]$.


Then, we can find a subgraph $Q$ in $N_6[v]$ defined in the proof of Lemma \ref{delete path}, \ref{eyeglass}, \ref{C4 2 deg 2}, \ref{triangle-free}, \ref{C4 contains deg 2}, \ref{C4 adjacent to deg 2}, \ref{C5 contains deg 2} or \ref{cubic} such that $G-Q$ does not contain a $K_{3,3}$ or $C_5$ as a component.
Put neighbors of $Q$ in $G$ into $T$, and remove $V(Q) \cap T$ from $T$.
Output the superlinear $L$-coloring of $G$ obtained from the superlinear coloring of $G-Q$ obtained by applying the subroutine with input $(G-Q,L,T)$ by extending it as in the proof of the corresponding lemma.
%

\noindent{\bf End of Description}

\bigskip

The correctness of the Subroutine 1 is clear.
Note that it is not hard to implement Subroutine 1 such that putting any vertex in $G$ into $T$ or removing any vertex from $T$ can be done in constant time.
To see the Subroutine 1 runs in time $O(\lvert V(G) \rvert)$, since $G$ has maximum degree at most three, it is sufficient to show that a superlinear $L$-coloring of a cycle $C$ with length other than five can be found in time $O(\lvert V(C) \rvert)$.
In fact, Skulrattanakulchai \cite{S06} proved that a $L$-coloring of $C^2$ can be found in $O(\lvert V(C) \rvert)$, and this coloring is a superlinear $L$-coloring of $C$.

Second, we introduce a subroutine to find a certain subgraph of a connected cubic graph.

\bigskip

\noindent{\bf Subroutine 2}

\noindent{\bf Input:} A connected cubic graph $G$ of order at least 11.

\noindent{\bf Output:} An induced subgraph $H$, where $H$ is isomorphic to $K_3$ or $K_{2,3}$, or $H$ is a cycle such that no vertex outside $H$ is adjacent to at least two vertices in $H$.

\noindent{\bf Running time:} $O(\lvert V(G) \rvert)$.

\noindent{\bf Description:}
Use brute force to find a triangle or a $K_{2,3}$ in $G$.
If we find a triangle or a $K_{2,3}$, then return it; otherwise, $G$ is triangle-free and $K_{2,3}$-free, and we pick a cycle $C$ in $G$. 

Let $C=v_0v_1...v_{k-1}v_0$ and let $u_i$ be the neighbor of $v_i$ other than $v_{i-1}$ and $v_{i+1}$, where the indices are computed modulo $k$.
For $i=0$ to $k-1$, consecutively see neighbors of $u_i$.
Whenever some $u_i$ is adjacent to at least two vertices in $C$, since $G$ is triangle-free and $K_{2,3}$-free, we can replace $C$ by $4$-cycle or a shorter cycle obtained from $u_i$ and the path in $C$ containing $v_{i+1}$ whose ends are in $N(u_i) \cap V(C)$.
So we can repeat this process until no $u_i$ is adjacent to at least two vertices in $C$.

\noindent{\bf End of Description}

\bigskip

The correctness of Subroutine 2 is clear.
We shall show that Subroutine 2 runs in time $O(\lvert V(G) \rvert)$.
Since $G$ has maximum degree three, it takes time $O(\lvert V(G) \rvert)$ to use brute force to find a triangle or a $K_{2,3}$.
Furthermore, finding a desired cycle mentioned in Subroutine 2 runs in linear-time since every vertex is visited by at most a bounded number of times.
So Subroutine 2 runs in linear-time.

Now, we give an algorithm to find a superlinear $L$-coloring of a connected subcubic graph $G$.

\bigskip

\noindent{\bf Algorithm for finding a superlinear $L$-coloring}

\noindent{\bf Input:} A connected subcubic graph $G$ which is not $K_{3,3}$ or $C_5$, and a family $L=\{L(v): v \in V(G)\}$ of lists of size $4$.

\noindent{\bf Output:} A superlinear $L$-coloring of $G$.

\noindent{\bf Running time:} $O(\lvert V(G) \rvert)$.

\noindent{\bf Description:}
If the order of $G$ is at most 10, then return a superlinear $L$-coloring of $G$ found by brute force.
So we assume that $G$ has order at least $11$.

Use breadth-first-search to check whether $G$ is cubic or not, and obtain the set $T$ consisting all vertices of degree at most two in $G$.
If $G$ is not cubic, then return the $L$-coloring obtained by applying Subroutine 1 with input $(G,L,T)$.
So we assume that $G$ is cubic.

Apply Subroutine 2 to obtain an induced subgraph $H$, where $H$ isomorphic to $K_3$ or $K_{2,3}$, or $H$ is a cycle such that no vertex outside $H$ is adjacent to at least two vertices in $H$.
Since $G$ is cubic, $G-H$ contains no $K_{3,3}$ or $C_5$ as a component if $H$ is isomorphic to $K_3$ or $K_{2,3}$.
If $G-H$ contains $C_5$ as a component, then $H$ is a cycle of length at least $6$ as $G$ has at least $11$ vertices, and we replace $H$ by this $5$-cycle in this case.
Hence, $G-H$ contains no $C_5$ or $K_{3,3}$ as a component.

If $H$ is a triangle or a $K_{2,3}$, respectively, then define $Q$ as in the proof of Lemma \ref{triangle-free} or \ref{K_{2,3}-free}, respectively; if $H$ is a cycle other than a triangle, then define $Q$ to be $H$.
Let $T'$ be the set of vetices of degree two in $G-Q$.
Note that $T'$ is a subset of $N(H)$, and $G-Q$ does not contain $C_5$ or $K_{3,3}$ as a component.
Return the superlinear $L$-coloring of $G$ obtained by extending the superlinear $L$-coloring obtained from applying Subroutine 1 with input $(G-Q, L, T')$ as the proof in Lemmas \ref{triangle-free}, \ref{K_{2,3}-free} and \ref{good cycle}.

%
\noindent{\bf End of the description.}

\bigskip

The correctness of the algorithm is clear.
Note that we only check whether a graph is cubic or not once, so the algorithm runs in linear-time as Subroutines 1 and 2 do.

\bigskip

\noindent{\bf Acknowledgement:}
The authors thank Professor Robin Thomas for many constructive suggestions and the referees for careful reading and valuable comments.

\end{document}